\documentclass[a4paper]{amsart}

\usepackage{amsmath, tabularx}
\usepackage{amsthm}
\usepackage{amssymb}
\usepackage{amsfonts}
\usepackage{paralist}
\usepackage{aliascnt}
\usepackage[initials]{amsrefs}
\usepackage{amscd}
\usepackage{blkarray}

\usepackage{setspace}
\usepackage[colorlinks=true,linkcolor=blue,urlcolor=blue]{hyperref}
\usepackage{tikz}
\usetikzlibrary{matrix}
\usetikzlibrary{arrows}

\newcommand{\RR}{\normalfont\mathbb{F}}
\newcommand{\kk}{\normalfont\mathbb{K}}

\newcommand{\ZZ}{\normalfont\mathbb{Z}}

\newcommand{\PP}{\normalfont\mathbb{P}}

\newcommand{\mm}{{\normalfont\mathfrak{m}}}

\newcommand{\QQ}{\widetilde{Q}}
\newcommand{\pp}{{\normalfont\mathfrak{p}}}

\newcommand{\HT}{\normalfont\text{ht}}

\newcommand{\HS}{{\normalfont\text{Hilb}}}

\newcommand{\Supp}{\normalfont\text{Supp}}
\newcommand{\Ass}{\normalfont\text{Ass}}
\newcommand{\Sym}{\normalfont\text{Sym}}
\newcommand{\Rees}{\mathcal{R}}

\newcommand{\Fitt}{\normalfont\text{Fitt}}
\newcommand{\EEQ}{\mathcal{K}}
\newcommand{\JJ}{\mathcal{J}}

\newcommand{\OO}{\mathcal{O}}

\newcommand{\LL}{\mathbb{L}}

\newcommand{\FF}{\normalfont\mathcal{F}}
\newcommand{\HL}{\normalfont\text{H}_{\mm}}
\newcommand{\HH}{\normalfont\text{H}}

\newcommand{\AAA}{\mathcal{A}}

\newcommand{\bideg}{\normalfont\text{bideg}}
\newcommand{\Proj}{\normalfont\text{Proj}}
\newcommand{\Spec}{\normalfont\text{Spec}}

\newtheorem{headthm}{Theorem}

\newaliascnt{headcor}{headthm}
\newtheorem{headcor}[headcor]{Corollary}
\aliascntresetthe{headcor}

\newaliascnt{corollary}{theorem}
\newtheorem{corollary}[corollary]{Corollary}
\aliascntresetthe{corollary}

\newaliascnt{lemma}{theorem}
\newtheorem{lemma}[lemma]{Lemma}
\aliascntresetthe{lemma}

\newaliascnt{conjecture}{theorem}

\aliascntresetthe{conjecture}

\newaliascnt{proposition}{theorem}
\newtheorem{proposition}[proposition]{Proposition}
\aliascntresetthe{proposition}

\newaliascnt{definition}{theorem}
\newtheorem{definition}[definition]{Definition}
\aliascntresetthe{definition}

\newaliascnt{notation}{theorem}
\newtheorem{notation}[notation]{Notation}
\aliascntresetthe{notation}

\newaliascnt{example}{theorem}

\aliascntresetthe{example}

\newaliascnt{examples}{theorem}

\aliascntresetthe{examples}

\newaliascnt{remark}{theorem}
\newtheorem{remark}[remark]{Remark}
\aliascntresetthe{remark}

\newaliascnt{problem}{theorem}

\aliascntresetthe{problem}

\newaliascnt{construction}{theorem}

\aliascntresetthe{construction}

\newaliascnt{algorithm}{theorem}

\aliascntresetthe{algorithm}

\newaliascnt{observation}{theorem}

\aliascntresetthe{observation}

\newaliascnt{defprop}{theorem}

\aliascntresetthe{defprop}

\def\equationautorefname~#1\null{(#1)\null}
\def\sectionautorefname~#1\null{Section #1\null}
\def\subsectionautorefname~#1\null{\S #1\null}

\begin{document}

\title[Saturated special fiber ring of height two perfect ideals]{Multiplicity of the saturated special fiber ring of height two perfect ideals}
\author{Yairon Cid-Ruiz}
\address{Department de Matem\`{a}tiques i Inform\`{a}tica, Facultat de Matem\`{a}tiques i Inform\`{a}tica, Universitat de Barcelona, Gran Via de les Corts Catalanes, 585; 08007 Barcelona, Spain.}
\email{ycid@ub.edu}
\urladdr{http://www.ub.edu/arcades/ycid.html}

\date{\today}

\subjclass[2010]{Primary 13A30; Secondary 14E05, 13D02, 13D45.}

\keywords{saturated special fiber ring, rational and birational maps, $j$-multiplicity, syzygies, Rees algebra, symmetric algebra, special fiber ring, multiplicity, Hilbert-Burch theorem, local cohomology.}

\thanks{The author was funded by the European Union's Horizon 2020 research and innovation programme under the Marie Sk\l{}odowska-Curie grant agreement No. 675789.
}

\begin{abstract}
	Let $R$ be a polynomial ring and $I \subset R$ be a perfect ideal of height two minimally generated by forms of the same degree.
	We provide a formula for the multiplicity of the \textit{saturated special fiber ring} of $I$. 
	Interestingly, this formula is equal to an elementary symmetric polynomial in terms of the degrees of the syzygies of $I$. 
 	Applying ideas introduced in \cite{MULTPROJ}, we obtain the value of the $j$-multiplicity of $I$ and an effective method for determining the degree and birationality of rational maps defined by homogeneous generators of $I$.
\end{abstract}
\maketitle
\vspace*{-.4cm}

\section{Introduction}

In \cite{MULTPROJ} an algebra called the \textit{saturated special fiber ring} was introduced, this algebra turns out to be an important tool in the study of rational and birational maps and is also related to the $j$-multiplicity of an ideal.
In this paper we compute the multiplicity of this algebra in the case of height two perfect ideals.
Interestingly, we express this multiplicity in terms of an elementary symmetric polynomial that depends on the degrees of the syzygies of the ideal.
As two simple corollaries, for this class of ideals, we obtain a closed formula for the $j$-multiplicity and an effective method for determining the degree and birationality of rational maps defined by homogeneous generators of these ideals.
 
\medskip

Let $\kk$ be a field, $R$ be the polynomial ring $R=\kk[x_0,x_1,\ldots x_r]$, and $\mm$ be the maximal irrelevant ideal $\mm=(x_0,x_1,\ldots,x_r)$.
Let $I\subset R$ be a perfect ideal of height two which is minimally generated by $s+1$ forms $\{f_0,f_1,\ldots,f_s\}$ of the same degree $d$.
As in \cite{MULTPROJ}, the saturated special fiber ring of $I$ is given by the algebra
$$
\QQ=\bigoplus_{n=0}^\infty {\left[\big(I^n:\mm^\infty\big)\right]}_{nd}.
$$
It can be seen as a saturated version of the classical special fiber ring.
To determine the multiplicity of $\QQ$, we need to study the first local cohomology module of the Rees algebra of $I$, and for this we  assume the condition $G_{r+1}$.
The condition $G_{r+1}$  means that $\mu(I_\pp) \le \dim(R_\pp)$ for every non-maximal ideal $\pp\in \Spec(R)$, where $\mu(I_\pp)$ denotes the minimal number of generators of $I_\pp$.
To study the Rees algebra one usually tries to reduce the problem in terms of the symmetric algebra, the assumption of $G_{r+1}$ is important in making possible this reduction.
After reducing the problem in terms of the symmetric algebra, we consider certain Koszul complex that provides an approximate resolution (see e.g. ~\cite{KPU_DEG_BOUND_LOC_COHOM}, \cite{CHARDIN_REG}) of the symmetric algebra,  and which permits us to compute the Hilbert series of $\QQ$.
By pursuing this general approach, we obtain the following theorem which is the main result of this paper.
   
\medskip

\begin{headthm}
		\label{MAIN_THM}
		Let $I\subset R=\kk[x_0,x_1,\ldots,x_r]$ be a homogeneous ideal minimally generated by $s+1$ forms $\{f_0,f_1,\ldots,f_s\}$ of the same degree $d$, where $s \ge r$.
		Suppose the following two conditions:
		\begin{enumerate}[(i)]
			\item $I$ is perfect of height two with Hilbert-Burch resolution of the form
			$$
			0 \rightarrow \bigoplus_{i=1}^sR(-d-\mu_i) \xrightarrow{\varphi} {R(-d)}^{s+1} \rightarrow I\rightarrow 0.			
			$$
			\item $I$ satisfies the condition $G_{r+1}$.
		\end{enumerate}
	Then, the multiplicity of the saturated special fiber ring $\QQ$ is given by
	$$
	e(\QQ)=e_r(\mu_1,\mu_2,\ldots,\mu_s),	
	$$
	where $e_r(\mu_1,\mu_2,\ldots,\mu_s)$ represents the $r$-th elementary symmetric polynomial
	$$
	e_r(\mu_1,\mu_2,\ldots,\mu_s) = \sum_{1\le j_1 < j_2 < \cdots < j_r \le s} \mu_{j_1}\mu_{j_2}\cdots\mu_{j_r}.
	$$
\end{headthm}

\medskip

As a first application of \autoref{MAIN_THM}, we obtain a closed formula for the $j$-multiplicity 
$$
j(I) = r
!\lim\limits_{n\rightarrow\infty} \frac{\dim_{\kk}\Big(\HH_{\mm}^0\left(I^n/I^{n+1}\right)\Big)}{n^r}
$$
of $I$.
The  $j$-multiplicity of an ideal was introduced in \cite{ACHILLES_MANARESI_J_MULT} and serves as a generalization of the Hilbert-Samuel multiplicity for non $\mm$-primary ideals.
It has applications in intersection theory (see \cite{FLENNER_O_CARROLL_VOGEL}), and the problem of finding formulas for it has been addressed in several papers (see e.g. ~\cite{JMULT_MONOMIAL, JEFFRIES_MONTANO_VARBARO, COMPUTING_J_MULT, POLINI_XIE_J_MULT}).
The following result gives a formula for the $j$-multiplicity of a whole family of ideals.

\medskip

\begin{headcor}
		\label{COR_J_MULT}
		Assume all the hypotheses and notations of \autoref{MAIN_THM}.
		Then, the $j$-multiplicity of $I$ is given by 
		$$
		j(I)=d\cdot e_r(\mu_1,\mu_2,\ldots,\mu_s).
		$$
\end{headcor}

\medskip

In the second application of \autoref{MAIN_THM}, we study the degree of a rational map $\FF:\PP^r \dashrightarrow \PP^s$ defined by the forms $f_0,f_1,\ldots,f_s$.
We show that the product of the degree of $\FF$ and the degree of the image of $\FF$ is equal to $e_r(\mu_1,\ldots,\mu_s)$.
From this we can determine the  degree of a rational map by just computing the degree of the image, and conversely, the degree of the map gives us the degree of the image.
In particular, we obtain that the map is birational if and only if the degree of the image is the maximum possible.
This idea of determining birationality by studying the syzygies of the base ideal is an active research topic (see e.g. ~\cite{AB_INITIO, Simis_cremona,KPU_blowup_fibers,EISENBUD_ULRICH_ROW_IDEALS,Hassanzadeh_Simis_Cremona_Sat,SIMIS_RUSSO_BIRAT,EFFECTIVE_BIGRAD,SIMIS_PAN_JONQUIERES,HASSANZADEH_SIMIS_DEGREES, HULEK_KATZ_SCHREYER_SYZ, MULTPROJ}).

\medskip

\begin{headcor}
		\label{COR_DEG_FF}
		Assume all the hypotheses and notations of \autoref{MAIN_THM}.
		Let $\FF$ be the rational map $\FF:\PP^r\dashrightarrow \PP^s$ given by 
		$$
		\left(x_0:\cdots:x_r\right) \mapsto 	\big(f_0(x_0,\ldots,x_r):\cdots:f_s(x_0,\ldots,x_r)\big),
		$$
		and $Y\subset \PP^s$ be the closure of the image of $\FF$.
		Then, the following two statements hold:
		\begin{enumerate}[(i)]
			\item $\deg(\FF)\cdot\deg_{\PP^s}(Y)=e_r(\mu_1,\mu_2,\ldots,\mu_s)$.
			\item $\FF$ is birational onto its image if and only if $\deg_{\PP^s}(Y)=e_r(\mu_1,\mu_2,\ldots,\mu_s)$.
		\end{enumerate}
\end{headcor}

\medskip

\section{Multiplicity of  the saturated special fiber ring}

The following notation will be assumed in the rest of this paper.

\begin{notation}
		\label{NOTA_SEC_2}
		Let $\kk$ be a field, $R$ be the polynomial ring $R=\kk[x_0,x_1,\ldots,x_r]$, and $\mm$ be the maximal irrelevant ideal $\mm=\left(x_0,x_1,\ldots,x_r\right)$. 
		Let $I$ be a homogeneous ideal minimally generated by $I=\left(f_0,f_1,\ldots,f_s\right) \subset R$ where $\deg(f_i)=d$ and $s\ge r$.
		Let $S$ be the polynomial ring $S=\kk[y_0,y_1,\ldots,y_s]$, and $\AAA$ be the bigraded polynomial ring $\AAA=R\otimes_{\kk}S=\kk[x_0,\ldots,x_r,y_0,\ldots,y_s]$.
		Let $Q$ be the standard graded $\kk$-algebra $Q=\kk[I_d]=\kk[f_0,f_1,\ldots,f_s]=\bigoplus_{n=0}^\infty {\left[I^n\right]}_{nd}$.
		
		We assume that $I$ is a perfect ideal of height two with Hilbert-Burch resolution of the form
		\begin{equation}
			\label{Hilb_Burch_presentation}
			0 \rightarrow \bigoplus_{i=1}^sR(-d-\mu_i) \xrightarrow{\varphi} {R(-d)}^{s+1} \rightarrow I\rightarrow 0.			
		\end{equation}
		We also suppose that $I$ satisfies the condition $G_{r+1}$, that is 
		$$
		\mu(I_{\pp}) \le \dim(R_{\pp}) \quad \text{ for all }\; \pp \in \Spec(R) \;\textit{ such that }\; \HT(\pp)<r+1.
		$$				
\end{notation}

\begin{remark}
	In terms of Fitting ideals, $I$ satisfies the condition $G_{r+1}$ if and only if $\HT(\Fitt_i(I)) > i$ for all $1 \le i < r+1$.
	So, from the presentation $\varphi$ of $I$, the condition $G_{r+1}$ is equivalent to $\HT(I_{r+1-i}(\varphi))>i$ for all $1 \le i < r+1$. 
	\begin{proof}
		It follows from \cite[Proposition 20.6]{EISEN_COMM}.
	\end{proof}
\end{remark}

We shall determine the multiplicity of the following algebra.

\begin{definition}[\cite{MULTPROJ}]
	The saturated special fiber ring of $I$ is given by the algebra
	$$
	\QQ=\bigoplus_{n=0}^\infty {\left[\big(I^n:\mm^\infty\big)\right]}_{nd}.
	$$
\end{definition}


The Rees algebra $\Rees(I)=\bigoplus_{n=0}^\infty I^nt^n \subset R[t]$ can be presented as a quotient of $\AAA$ by using the map
\begin{eqnarray*}
\label{presentation_Rees}
\Psi: \AAA & \longrightarrow & \Rees(I) \subset R[t] \\ \nonumber
y_i & \mapsto & f_it.
\end{eqnarray*}
We set $\bideg(x_i)=(1,0)$, $\bideg(y_j)=(0,1)$ and $\bideg(t)=(-d, 1)$, which implies that $\Psi$ is bihomogeneous of degree zero, and so $\Rees(I)$ has a structure of bigraded $\AAA$-algebra.
If $M$ is a bigraded $\AAA$-module and $c$ a fixed integer, then we write
$$
{\left[M\right]}_c = \bigoplus_{n  \in \ZZ} M_{(c,n)}.
$$
We remark that ${\left[M\right]}_c$ has a natural structure as a graded $S$-module.

As noted in \cite{MULTPROJ}, to study the algebra $\QQ$ it is enough to consider the degree zero part in the $R$-grading of the bigraded $\AAA$-module $\HL^1\left(\Rees(I\right))$ (see e.g. ~\cite[Lemma 2.1]{DMOD}).

\begin{remark}
	Let $X$ be the scheme $X=\Proj_{R\text{-gr}}\left(\Rees(I)\right)$, where $\Rees(I)$ is only considered as a graded $R$-algebra.
	From \cite[Theorem A4.1]{EISEN_COMM}, we obtain the following short exact sequence
	$$
	0 \rightarrow {\left[\Rees(I)\right]}_0 \rightarrow \HH^0(X, \OO_X) \rightarrow {\left[\HL^1\left(\Rees(I)\right)\right]}_0 \rightarrow 0.
	$$
	By identifying $Q\cong {\left[\Rees(I)\right]}_0$ and $\QQ \cong \HH^0(X, \OO_X)$, we obtain the short exact sequence
	\begin{equation}
		\label{EQ_SPECIAL_AND_SAT}
		0 \rightarrow Q \rightarrow \QQ \rightarrow {\left[\HL^1\left(\Rees(I)\right)\right]}_0 \rightarrow 0.
	\end{equation}
\end{remark}

\begin{remark}
	\label{REM_f.g. Q-mod}
	From \cite[Proposition 2.7$(i)$, Lemma 2.8$(ii)$]{MULTPROJ} we have that $\QQ$ and ${\left[\HL^1\left(\Rees(I)\right)\right]}_0$ have natural structures of finitely generated $Q$-modules.
\end{remark}

The Rees algebra is a very difficult object to study, but, under the present conditions, we have that the module ${\left[\HL^1\left(\Rees(I)\right)\right]}_0$ coincides with ${\left[\HL^1\left(\Sym(I)\right)\right]}_0$ (see \autoref{LemSymProps}$(iii)$).
So, the main idea is to bypass the Rees algebra and consider the symmetric algebra instead.
From the presentation \autoref{Hilb_Burch_presentation} of $I$, we obtain the ideal 
$$
\JJ = \left(g_1,\ldots,g_s\right) = I_1\big(\left[y_0,\ldots,y_s\right]\cdot\varphi\big)
$$
of defining equations of the symmetric algebra.
Thus, $\Sym(I)$ is a bigraded $\AAA$-algebra presented by the quotient
$$
\Sym(I) \cong \AAA/\JJ.
$$
We have the following canonical short exact sequence relating both algebras
\begin{equation}
\label{EqSymRees}
0 \rightarrow \EEQ \rightarrow \Sym(I) \rightarrow \Rees(I) \rightarrow 0,	
\end{equation}
where $\EEQ$ is the $R$-torsion submodule of $\Sym(I)$.

We will consider the Koszul complex $\LL_\bullet=K_\bullet(g_1,\ldots,g_s;\AAA)$ associated to $\{g_1,\ldots,g_s\}$:
$$
\LL_\bullet: \quad 0 \rightarrow \LL_s \rightarrow \cdots \rightarrow \LL_i \rightarrow \cdots \rightarrow \LL_1 \rightarrow \LL_0 
$$
where 
\begin{equation}
	\label{EQ_KOSZUL_COMP_SYM}
	\LL_i = \bigwedge\nolimits^{\!i}\left(\bigoplus_{j=1}^s\AAA(-\mu_j, -1)\right).	
\end{equation}
This complex will not be exact in general, but the homology modules will have small enough Krull dimension.
It will give us an ``approximate resolution'' of the symmetric algebra (see e.g. ~\cite{KPU_DEG_BOUND_LOC_COHOM}, \cite{CHARDIN_REG}), from which we can read everything we need.

In the following lemma we gather some well-known properties of $\Sym(I)$ under the present conditions,  we include them for the sake of completeness.

\begin{lemma}
	\label{LemSymProps}
	Using \autoref{NOTA_SEC_2}, the following statements hold:
	\begin{enumerate}[(i)]
		\item $\dim\left(\Sym(I)\right)=\max\left(\dim(R)+1, \mu(I)\right)=\max\left(r+2, s+1\right)$.
		\item $\EEQ=\HL^0\left(\Sym(I)\right)$.
		\item $\HL^i\left(\Rees(I)\right) \cong \HL^i\left(\Sym(I)\right)$ for all $i\ge 1$.
		\item If $s \le r +1$, then $\Sym(I)$ is a complete intersection.
		\item For all $s\ge1$,  $\Sym(I)$ is a complete intersection on the punctured spectrum of $R$.
	\end{enumerate}
	\begin{proof}
		$(i)$ Follows from the dimension formula for symmetric algebras (see \cite{HUNEKE_ROSSI}, \cite[Theorem 1.2.1]{vasconcelos1994arithmetic}) and the condition $G_{r+1}$.
		
		$(ii)$ It follows from \cite[Corollary 4.8 , \S 5]{HSV_TRENTO_SCHOOL} (also, see \cite[\S 3.7]{KPU_GOR3})).
		
		$(iii)$ For each $i \ge 1$, the short exact sequence \autoref{EqSymRees} yields the long exact sequence 
		$$
		\HL^i\left(\EEQ\right) \rightarrow \HL^i\left(\Sym(I)\right) \rightarrow \HL^i\left(\Rees(I)\right) \rightarrow \HL^{i+1}\left(\EEQ\right).
		$$
		From part $(ii)$ and \cite[Corollary 2.1.7]{Brodmann_Sharp_local_cohom}, we have that $\HL^i\left(\EEQ\right)=\HL^{i+1}\left(\EEQ\right)=0$, and so we obtain the required isomorphism.
		
		$(iv)$ Using part $(i)$, in this case we have that $\dim(\Sym(I))=r+2$.
		Hence, we get 
		$$
		\HT(\JJ)=\dim(\AAA) - (r+2)=(r+s+2)-(r+2)=s=\mu(\JJ),		
		$$
		and so $\Sym(I)$ is a complete intersection.
		
		$(v)$ For each $\pp \in \Spec(R)$ such that $\HT(\pp)<r+1$, the same argument of part $(i)$ now yields that $\dim({\Sym(I)}_\pp)=\dim(R_\pp)+1$.
		Thus, we have
		$$
		\HT(\JJ_\pp) = \dim(\AAA_\pp) - \dim({\Sym(I)}_\pp)= \dim(R_\pp) + s+1 - (\dim(R_\pp)+1) = s=\mu(\JJ_\pp).
		$$
		Then, for $i\ge 1$, the homology module $\HH_i(\LL_\bullet)$ is supported on the maximal ideals of $\Spec(R)$, but since the associated primes $\Ass_R(\HH_i(\LL_\bullet))$ are homogeneous, it necessarily gives that $\Supp_R(\HH_i(\LL_\bullet))=\{\mm\}$.
		Therefore, ${\Sym(I)}_\pp$ is a complete intersection for $\pp \in \Spec(R) \setminus \{\mm\}$.
	\end{proof}
\end{lemma}

The restriction to degree zero part in the $R$-grading of the equality $\EEQ=\HL^0\left(\Sym(I)\right)$ (\autoref{LemSymProps}$(ii)$) and the short exact sequence \autoref{EqSymRees} yield the following
\begin{equation}
	\label{EQ_SymRees_deg0}
	0 \rightarrow {\left[\HL^0\left(\Sym(I)\right)\right]}_0 \rightarrow S \rightarrow Q \rightarrow 0,
\end{equation}
under the identifications ${\left[\Sym(I)\right]}_0=S$ and ${\left[\Rees(I)\right]}_0=Q$.

The next proposition will be an important technical tool.

\begin{proposition}
	\label{PropLocCohoSym}
	Assume \autoref{NOTA_SEC_2}. Then, we have the following isomorphisms of bigraded $\AAA$-modules
	$$
	\HH_{i}\Big(\HL^{r+1}(\LL_\bullet)\Big) \cong \begin{cases}
		\HL^{r+1-i}\left(\Sym(I)\right) \quad \text{ if } i \le r+1\\
		\HH_{i-r-1}(\LL_\bullet) \qquad\quad\;\, \text{ if } i \ge r+2,
	\end{cases}	
	$$		
	where $\HL^{r+1}(\LL_\bullet)$ represents the complex obtained after applying the functor $\HL^{r+1}(\bullet)$ to $\LL_\bullet$.
	\begin{proof}
		Let $\mathbb{G}^{\bullet,\bullet}$ be the first quadrant double complex given by $\mathbb{G}^{p,q}=\LL_{s-p}\otimes_R C_\mm^q$, where $C_\mm^\bullet$ is the \v{C}ech complex corresponding with the maximal irrelevant ideal $\mm$.
		
		Since we have that
		\begin{equation}
			\label{EQ_LOCAL_COHOM}
			\HL^{p}(\AAA) \cong \begin{cases}
			\frac{1}{x_0x_1\cdots x_r}\kk[x_0^{-1},x_1^{-1},\ldots,x_r^{-1}] \otimes_{\kk} S \quad \text{ if } p = r+1\\
			0 \;\;\;\quad\qquad\qquad\qquad\qquad\qquad\qquad\quad\;\;\; \text{otherwise},
			\end{cases}			
		\end{equation}
		then the spectral sequence coming from the first filtration is given by 
		$$
		{}^{\text{I}}E_1^{p,q}=\begin{cases}
			\HL^{r+1}(\LL_{s-p}) \quad \text{ if } q = r+1\\
			0 \quad\qquad\qquad\;\; \text{ otherwise}.
		\end{cases}
		$$
		
		On the other hand, \autoref{LemSymProps}$(v)$ implies that ${\left(\LL_\bullet\right)}_\pp$ is exact for all $\pp \in \Spec(R) \setminus \{\mm\}$.
		So, for all $i \le s-1$, $\HH_{s-i}(\LL_\bullet)$ is supported on $V(\mm)$ and the Grothendieck vanishing theorem (see e.g. ~\cite[Theorem 6.1.2]{Brodmann_Sharp_local_cohom})  implies that 
		$$
		\HL^j(\HH_{s-i}(\LL_\bullet))=0
		$$ 
		for all $j\ge 1$.
		Also, we have that 
		$$
		\HL^0(\HH_{s-i}(\LL_\bullet)) = \HH_{s-i}(\LL_\bullet)
		$$
		for $i\le s-1$.
		Therefore, the spectral sequence corresponding with the second filtration is given by
		$$
		{}^{\text{II}}E_2^{p,q} \cong \begin{cases}
				\HL^p\left(\Sym(I)\right) \quad\quad \text{ if } q=s \\
				\HH_{s-q}(\LL_\bullet) \qquad\;\,\quad\, \text{ if } p = 0 \text{ and } q \le s-1\\
				0 \quad\qquad\qquad\qquad\; \text{ otherwise}.
		\end{cases}
		$$
		
		Finally, from the convergence of both spectral sequences we obtain the following isomorphisms of bigraded $\AAA$-modules
		$$
			\HH_{i}\left(\HL^{r+1}(\LL_\bullet)\right) \cong\HH^{r+1+s-i}\left(\text{Tot}(\mathbb{G}^{\bullet,\bullet})\right) 
			\cong \begin{cases}
			\HL^{r+1-i}\left(\Sym(I)\right) \quad \text{ if } i \le r+1\\
			\HH_{i-r-1}(\LL_\bullet) \qquad\quad\;\, \text{ if } i \ge r+2
			\end{cases}
		$$
		for all $i \ge 0$.		
	\end{proof}
\end{proposition}

The following lemma contains some dimension computations that will be needed in the proof of \autoref{MAIN_THM}.
The first one shows that $I$ has maximal analytic spread and it is obtained directly from  \cite{Ulrich_Vasc_Eq_Rees_Lin_Present}.
The second one is a curious interplay between the algebraic properties of $I$ and the geometric features of the corresponding rational map \autoref{EQ_RAT_MAP}, that follows from \cite[Proposition 3.1]{MULTPROJ}.
(A much stronger generalization of \cite[Proposition 3.1]{MULTPROJ} was recently obtained in \cite[Theorem 4.4]{DEGREE_SPECIALIZATION}.)

\begin{lemma}
		\label{LEM_DIM_LOC_SYM_TWO}
		Using \autoref{NOTA_SEC_2}, the following statements hold:
		\begin{enumerate}[(i)]
			\item $\ell(I)=\dim\left(\Rees(I)/\mm \Rees(I)\right)=r+1$.
			\item $\dim(Q)=r+1$.
			\item The corresponding rational $\FF:\PP^r \dashrightarrow \PP^s$ in \autoref{COR_DEG_FF} is generically finite.
			\item $
			\dim\left({\left[\HL^i\left(\Sym(I)\right)\right]}_0\right) \,\le\, r$\; for all $i \ge 2$.
		\end{enumerate}

		\begin{proof}
			$(i)$
			In the case $r=s$, we get from \cite[Theorem 4.1]{Ulrich_Vasc_Eq_Rees_Lin_Present} that $I$ is of linear type and so $\dim\left(\Rees(I)/\mm \Rees(I)\right)=r+1$.
			When $s\ge r+1$, then the result follows from \cite[Corollary 4.3]{Ulrich_Vasc_Eq_Rees_Lin_Present}.
			
			$(ii)$ Since we have $$
			\Rees(I) = {\left[\Rees(I)\right]}_0\; \bigoplus\;\left( \bigoplus_{n=1}^\infty{\left[\Rees(I)\right]}_n\right)=Q \;\bigoplus\; \mm \Rees(I), 
			$$
			then we get an isomorphism $Q \cong \Rees(I)/\mm \Rees(I)$ of graded $\kk$-algebras.
			Thus, from part $(i)$, $\dim(Q)=\dim\left(\Rees(I)/\mm\Rees(I)\right)=r+1$.
			
			$(iii)$ Let $Y$ be the closure of the image of $\FF$. 
			Since $Q=\kk[f_0,\ldots,f_s]$ corresponds with the homogeneous coordinate ring of $Y$, the claim follows from part $(ii)$ and \cite[Corollary 3.3, Proposition 3.14]{DEGREE_SPECIALIZATION}.
			
			$(iv)$
			Let $i\ge 2$.
			From \autoref{LemSymProps}$(iii)$, we have ${\left[\HL^i\left(\Sym(I)\right)\right]}_0 \cong {\left[\HL^i\left(\Rees(I)\right)\right]}_0$.
			The rational map $\FF$ is generically finite due to part $(iii)$, and so the inequality follows directly from \cite[Proposition 3.1]{MULTPROJ}.
		\end{proof}
\end{lemma}

Now we are ready for the proof of the main theorem.

\begin{proof}[Proof of \autoref{MAIN_THM}]
		The whole point of this proof is to analyze the homology modules of the complex
		$$
		\RR_\bullet={\left[\HL^{r+1}(\LL_\bullet)\right]}_0\,: \quad 0 \rightarrow {\left[\HL^{r+1}(\LL_s)\right]}_0 \rightarrow \cdots \rightarrow {\left[\HL^{r+1}(\LL_1)\right]}_0 \rightarrow {\left[\HL^{r+1}(\LL_0)\right]}_0
		$$
		obtained by applying  $\HL^{r+1}(\bullet)$ to the complex $\LL_\bullet$ and then restricting to the degree zero part in the $R$-grading.
		From \autoref{EQ_KOSZUL_COMP_SYM} and \autoref{EQ_LOCAL_COHOM}, we can make the identification 
		$$
		\RR_i = {\left[\HL^{r+1}(\LL_i)\right]}_0 \,\cong\, {S(-i)}^{m_i},		
		$$
		where
		$$
		m_i = \sum_{1\le j_1 < \cdots < j_i \le s} \binom{\sum_{e=1}^i\mu_{j_e} - 1}{r}.
 		$$
 		
 		First, from \autoref{PropLocCohoSym} we have 
 		\begin{equation*} 
 			 \HH_i(\RR_\bullet) \,\cong\, {\left[\HH_{i-r-1}(\LL_\bullet)\right]}_0 \;\text{ for }\; i \ge r+2,
 		\end{equation*}
 		then the fact that ${\left[\LL_k\right]}_0=0$ for $k\ge 1$ (see \autoref{EQ_KOSZUL_COMP_SYM}) yields the vanishing
		\begin{equation}
			\label{EQ_VANISH_GE_r+2}
			 \HH_i(\RR_\bullet)= 0 \;\text{ for all }\; i \ge r+2.
		\end{equation}
 		On the other hand,  \autoref{PropLocCohoSym} also gives that 
 		$$
 		\HH_{i}(\RR_\bullet) \cong {\left[\HL^{r+1-i}(\Sym(I))\right]}_0 \quad \text{ for } i \le r+1,
 		$$
 		and \autoref{LEM_DIM_LOC_SYM_TWO}$(iv)$ implies that
 		\begin{equation}
 			\label{EQ_DIM_LE_r-1}
			\dim\left(\HH_i(\RR_\bullet)\right)\le r \;\text{ for all } i\le r-1.
 		\end{equation} 
 		
 		Let $B_\bullet$, $Z_\bullet$ and $H_\bullet$ be the boundaries, cycles and homology modules of the complex $\RR_\bullet$, respectively.
		We have the following short exact sequences 
		\begin{align*}
			0 \rightarrow B_i \rightarrow &Z_i \rightarrow H_i \rightarrow 0\\
			0 \rightarrow Z_i \rightarrow &\RR_i \rightarrow B_{i-1} \rightarrow 0 
		\end{align*}
		for all $i$.
		By using the additivity of Hilbert series and assembling all these short exact  sequences we obtain the following equation
		$$
		\sum_{i=0}^s{(-1)}^i\HS_{H_i}(T) = \sum_{i=0}^s {(-1)}^i\HS_{\RR_i}(T) .
		$$
		Using \autoref{EQ_VANISH_GE_r+2} and \autoref{EQ_DIM_LE_r-1}, it follows that $\HS_{H_i}(T)=0$ for $i\ge r+2$, and that we can write 
		$$
		\HS_{H_i}(T)=\frac{G_i(T)}{{(1-T)}^{e_i}} \quad\text{ for } i \le r-1
		$$
		where $G_i(T) \in \ZZ[T]$ and $e_i=\dim(H_i)\le r$ (see e.g. ~\cite[Section 4.1]{BRUNS_HERZOG}).
		Therefore, we obtain the following equation
		$$
		\frac{C(T)}{{(1-T)}^{s+1}}  + {(-1)}^{r}\HS_{H_r}(T) + {(-1)}^{r+1}\HS_{H_{r+1}}(T)= \frac{G(T)}{{(1-T)}^{s+1}}
		$$
		where 
		$$
		C(T)=\sum_{i=0}^{r-1}{(-1)}^i{(1-T)}^{s+1-e_i}G_i(T) \;\;\text{ and }\;\; G(T)=\sum_{i=0}^s{(-1)}^im_iT^i.
		$$
		The isomorphisms of \autoref{PropLocCohoSym} yield that 
		\begin{equation}
			\label{EQ_Hm1Sym_Hilb}
			\HS_{{\left[\HL^1\left(\Sym(I)\right)\right]}_0}(T) = \HS_{{\left[\HL^0\left(\Sym(I)\right)\right]}_0}(T) + \frac{{(-1)}^{r}G(T)+{(-1)}^{r+1}C(T)}{{(1-T)}^{s+1}}
		\end{equation}
		From the short exact sequence \autoref{EQ_SymRees_deg0} we obtain that
		\begin{equation}
			\label{EQ_Hm0Sym_Hilb}
			\HS_{{\left[\HL^0\left(\Sym(I)\right)\right]}_0}(T) = \HS_S(T)-\HS_Q(T)=\frac{1}{{(1-T)}^{s+1}}-\HS_Q(T), 
		\end{equation}
		and the short exact sequence \autoref{EQ_SPECIAL_AND_SAT} and \autoref{LemSymProps}$(iii)$ yield that
		\begin{equation}
			\label{EQ_QQ_Hilb}
			\HS_{\QQ}(T) = \HS_{Q}(T) + \HS_{{\left[\HL^1\left(\Sym(I)\right)\right]}_0}(T).
		\end{equation}	
		Hence, by summing up \autoref{EQ_Hm1Sym_Hilb}, \autoref{EQ_Hm0Sym_Hilb} and \autoref{EQ_QQ_Hilb} we get 
		$$
		\HS_{\QQ}(T) = \frac{1+{(-1)}^{r}G(T)+{(-1)}^{r+1}C(T)}{{(1-T)}^{s+1}}.
		$$		
		
		Let $F(T)=1+{(-1)}^{r}G(T)+{(-1)}^{r+1}C(T)$.
		Since $Q \hookrightarrow \QQ$ is an integral extension (see \autoref{REM_f.g. Q-mod}), it follows that $\dim(\QQ)=\dim(Q)$.
		From \autoref{LEM_DIM_LOC_SYM_TWO}$(ii)$ we have that $\dim(\QQ)=\dim(Q)=r+1$, then well-known properties of Hilbert series (see e.g. ~\cite[Section 4.1]{BRUNS_HERZOG}) give us that
		$$
		F(T)={(1-T)}^{s-r}F_1(T),
		$$
		where $F_1(1)\neq 0$ and $e(\QQ)=F_1(1)$.
		The fact that $e_i\le r$ for $i\le r-1$, implies that $C^{(s-r)}(1)=0$.
		By denoting 
		$$
		P(T)=1+{(-1)}^rG(T)= 1+\sum_{i=0}^s{(-1)}^{r+i}m_iT^i,
		$$
		we get $P^{(s-r)}(1)=F^{(s-r)}(1)$, and so by taking the $(s-r)$-th derivatives of $F(T)$ and $P(T)$ we obtain that
		\begin{align*}
			{(-1)}^{s-r}(s-r)!\cdot F_1(1) &= P^{(s-r)}(1)\\
			&=\begin{cases}
				  1+\sum_{i=0}^r{(-1)}^{r+i}m_i \;\;\;\;\;\qquad\qquad\text{ if } s = r\\
				  \sum_{i=s-r}^s{(-1)}^{r+i}m_i(s-r)!\binom{i}{s-r} \;\; \text{ if } s > r.
				 \end{cases}
		\end{align*}
		The substitution of $e(\QQ)=F_1(1)$ gives us that
		\begin{equation}
			\label{EQ_value_of_mult}
			e(\QQ)=\begin{cases}
			1+\sum_{i=0}^r{(-1)}^{r+i}m_i \qquad\qquad\; \text{ if } s = r\\
			\sum_{i=s-r}^s{(-1)}^{s+i}m_i\binom{i}{s-r} \qquad\;\;\, \text{ if } s > r.
			\end{cases}
		\end{equation}
		Finally, the formula of the theorem is obtained from \autoref{LEM_formulas}$(iii),(iv)$ below.
\end{proof}

In the following lemma we use simple combinatorial techniques to reduce the equation \autoref{EQ_value_of_mult}.

\begin{lemma}
		\label{LEM_formulas}
		The following formulas hold:
		\begin{enumerate}[(i)]
			 \item For $0 \le k \le r$,
			 $$
			 \sum_{i=\max\{k,s-r\}}^{s}{(-1)}^i\binom{i}{s-r}\binom{s-k}{i-k} = \begin{cases}
					{(-1)}^s \quad\text{ if } k = r\\
					0 \qquad\;\;\;\,\text{ if } k < r.
			 \end{cases}
			 $$			  
			 \item For $1 \le \ell \le r$, 
			 $$
			 \sum_{i=s-r}^s{(-1)}^i\binom{i}{s-r}\sum_{1\le j_1 < \cdots < j_i \le s}\Big(\sum_{e=1}^i\mu_{j_e}\Big)^\ell = \begin{cases}
				{(-1)}^sr!\cdot e_r(\mu_1,\ldots,\mu_s)  \;\;\text{ if } \ell = r\\
				0 \qquad\qquad\qquad\qquad\quad\;\;\;\text{ if } \ell < r.
			 \end{cases}
			 $$
			 \item For $s>r$,			 
			 $$
			 \sum_{i=s-r}^s{(-1)}^i\binom{i}{s-r}\sum_{1\le j_1 < \cdots < j_i \le s}\binom{\sum_{e=1}^i\mu_{j_e}-1}{r}={(-1)}^s\cdot e_r(\mu_1,\ldots,\mu_s).
			 $$
			 \item For $s=r$,
			 $$
			 1+\sum_{i=0}^r{(-1)}^{i+r}\sum_{1\le j_1 < \cdots < j_i \le r}\binom{\sum_{e=1}^i\mu_{j_e}-1}{r} = \mu_1\mu_2\cdots\mu_r.
			 $$
		\end{enumerate}
		\begin{proof}
			$(i)$ We depart from the identity 
			$$
			(1-T)^{s-k}T^k= \sum_{i=k}^{s}{(-1)}^{i-k}\binom{s-k}{i-k}T^i,
			$$
			then by taking the $(s-r)$-th derivative in both sides we get 
			$$
			 {\Big((1-T)^{s-k}T^k\Big)}^{(s-r)}= \sum_{i=\max\{k,s-r\}}^s{(-1)}^{i-k}\binom{s-k}{i-k}(s-r)!\binom{i}{s-r}T^{i-s+r}.
			$$
			Since $s-k\ge s-r$, the substitution $T=1$ yields the result. 
			
			$(ii)$ For each set of indexes $\{j_1,\ldots,j_i\}$ we have 
			\begin{equation}
				\label{EQ_MULTINOMIAL}
				\Big(\sum_{e=1}^i\mu_{j_e}\Big)^\ell=\sum_{\ell_1+\cdots+\ell_i=\ell} \binom{\ell}{\ell_1,\ldots,\ell_i}\mu_{j_1}^{\ell_1}\cdots \mu_{j_i}^{\ell_i}.
			\end{equation}
			We will proceed by determining the coefficients of each of the monomials $\mu_{j_1}^{\ell_1}\cdots \mu_{j_i}^{\ell_i}$ in the equation. 
			Since $\binom{\ell}{\ell_1,\ldots,\ell_i}=\binom{\ell}{\ell_1,\ldots,\ell_i,0}$, we can consider the case where $\ell_1\neq 0,\ldots,\ell_k\neq0$.
			
			We fix $1 \le k \le r$ and the monomial $\mu_{i_1}^{b_1}\cdots\mu_{i_k}^{b_k}$ where $b_1\neq 0,\ldots,b_k\neq0$ and $b_1+\cdots+b_k=\ell$.
			For each set of indexes $\{j_1,\ldots,j_i\} \supset \{i_1,\ldots,i_k\}$, the monomial $\mu_{i_1}^{b_1}\cdots\mu_{i_k}^{b_k}$ appears once in the equation \autoref{EQ_MULTINOMIAL}, and the number of these sets is equal to $\binom{s-k}{i-k}$.
			Thus, for each $i \ge k$, the coefficient of $\mu_{i_1}^{b_1}\cdots\mu_{i_k}^{b_k}$ in the expression 
			$$
			\sum_{1\le j_1 < \cdots < j_i \le s}\Big(\sum_{e=1}^i\mu_{j_e}\Big)^\ell
			$$
			is equal to $\binom{s-k}{i-k}\binom{\ell}{b_1,\ldots,b_k}$.
			So the total coefficient of $\mu_{i_1}^{b_1}\cdots\mu_{i_k}^{b_k}$ is given by 
			$$
			\binom{\ell}{b_1,\ldots,b_k}\sum_{i=\max\{k,s-r\}}^{s}{(-1)}^i\binom{i}{s-r}\binom{s-k}{i-k}.
			$$
			From part $(i)$, we have that this coefficient vanishes when $k<r$ and that it is equal to ${(-1)}^sr!$ when $k = r$ because $\ell \le r$.
		
			Therefore, for $\ell<r$ we have that the equation vanishes, and for $\ell=r$ that the only monomials in the equation are those of the elementary symmetric polynomial $e_r(\mu_1,\ldots,\mu_s)$ and the coefficient of all of them is ${(-1)}^sr!$.
			
			$(iii)$ We can write 
			\begin{align}
				\label{EQ_BINOMIAL_MU_S}
				\binom{\sum_{e=1}^i\mu_{j_e}-1}{r}&=\frac{\left(\sum_{e=1}^i\mu_{j_e}-1\right)\left(\sum_{e=1}^i\mu_{j_e}-2\right)\cdots\left(\sum_{e=1}^i\mu_{j_e}-r\right)}{r!}\\
				\nonumber
				&=\frac{1}{r!}\sum_{\ell=0}^r{(-1)}^{r-\ell}e_{r-\ell}(1,2,\ldots,r){\left(\sum_{e=1}^i\mu_{j_e}\right)}^\ell.		
			\end{align}
			Therefore, by summing up and using part $(ii)$, we obtain the required formula.
			
			$(iv)$ From  equation \autoref{EQ_BINOMIAL_MU_S} and part $(ii)$ we have
			$$
			\sum_{i=0}^r{(-1)}^{i+r}\sum_{1\le j_1 < \cdots < j_i \le r}\binom{\sum_{e=1}^i\mu_{j_e}-1}{r}=\mu_1\mu_2\cdots \mu_r + \sum_{i=1}^r{(-1)}^i\binom{r}{i}.
			$$
			Thus we get the result from the identity $\sum_{i=0}^r{(-1)}^i\binom{r}{i}=0$.
		\end{proof}
\end{lemma}

From the main theorem we easily obtain a closed formula for the $j$-multiplicity of $I$.
\begin{proof}[Proof of \autoref{COR_J_MULT}]
	From \cite[Lemma 2.10]{MULTPROJ} we have that $j(I)=d\cdot e(\QQ)$, then the result follows from the computation of \autoref{MAIN_THM}.
\end{proof}

\section{Degree of rational maps}

In this short section we study the degree of the rational map 
\begin{align}
\label{EQ_RAT_MAP}
&\FF : \PP^r\dashrightarrow \PP^s\\
\nonumber
\left(x_0:\cdots:x_r\right) &\mapsto 	\big(f_0(x_0,\ldots,x_r):\cdots:f_s(x_0,\ldots,x_r)\big),
\end{align}
whose base ideal $I=(f_0,f_1,\ldots,f_s)$ satisfies all the conditions of \autoref{NOTA_SEC_2}.
Here we obtain a suitable generalization of \cite[Theorem 4.9 $(1),(2)$]{KPU_blowup_fibers}, where we relate the degree of $\FF$ and the degree of its image with the formula obtained in \autoref{MAIN_THM}.
An interesting result is that $\FF$ is birational onto its image if and only if the degree of the image is the maximum possible.

Let $Y\subset \PP^s$ be the closure of the image of $\FF$.
From \autoref{LEM_DIM_LOC_SYM_TWO}$(iii)$ we have that $\FF$ is generically finite, and that the degree of $\FF$ is equal to the dimension of the field extension 
$$
\deg(\FF) = \left[K(\PP^r):K(Y)\right],
$$
where $K(\PP^r)$ and $K(Y)$ represent the fields of rational functions of $\PP^r$ and $Y$, respectively.

The main result of this section is a simple corollary of \cite{MULTPROJ} and \autoref{MAIN_THM}.

\begin{proof}[Proof of \autoref{COR_DEG_FF}]
	 From \cite[Theorem 2.4$(iii)$]{MULTPROJ} we have that $e(\QQ)=\deg(\FF)\cdot \deg_{\PP^s}(Y)$, then the result is obtained from the computation of \autoref{MAIN_THM}.
\end{proof}

We have that in the literature special cases of \autoref{COR_DEG_FF} have appeared before.
For instance, in \cite[Proposition 5.3]{COX_EQ_PARAM} a particular case of \autoref{COR_DEG_FF} was obtained for parameterized surfaces.
In the following simple corollaries, we prove the same result of \cite[Theorem 4.9 $(1),(2)$]{KPU_blowup_fibers}, and we generalize \cite[Proposition 5.2]{MULTPROJ}.

\begin{corollary}
	With the same notations above, if $r=1$, i.e. $\FF$ is of the form $\FF:\PP^1\dashrightarrow \PP^s$, then 
	$
	\deg(\FF)\cdot\deg_{\PP^s}(Y)=d.
	$
	\begin{proof}
		From the Hilbert-Burch theorem (see e.g. ~\cite[Theorem 20.15]{EISEN_COMM}), in \autoref{NOTA_SEC_2},  $I$ is minimally generated by the maximal minors of $\varphi$.
		Therefore, we have that $d=\mu_1+\mu_2+\cdots+\mu_s=e_1(\mu_1,\mu_2,\ldots,\mu_s)$.
	\end{proof}
\end{corollary}

\begin{corollary}
		With the same notations above, if $r=s$, i.e. $\FF$ is of the form $\FF:\PP^r\dashrightarrow\PP^r$, then $\deg(\FF)=\mu_1\mu_2\cdots\mu_r.$
		\begin{proof}
			In this case we have $Y=\PP^r$ and so $\deg_{\PP^r}(Y)=1$.
			Hence the equality follows from the fact that $e_r(\mu_1,\mu_2,\ldots,\mu_r)=\mu_1\mu_2\cdots \mu_r.
			$
		\end{proof}
\end{corollary}

\section*{Acknowledgments}
\sloppy
The author thanks Frank-Olaf Schreyer for helpful discussions during the Macaulay2 workshop in Leipzig. 
The author is thankful to Laurent Bus\'e and Carlos D'Andrea for useful suggestions and discussions.
The author thanks Laurent Bus\'e for pointing out the previous result of \cite[Proposition 5.3]{COX_EQ_PARAM}.
The use of \textit{Macaulay2} \cite{MACAULAY2} was very important in the preparation of this paper.
The author wishes to thank the referee for several suggestions to improve the exposition.


\end{document}